\documentclass[11pt]{imsart}

\usepackage{amsfonts}
\usepackage[dvips]{graphics}
\usepackage{amsmath}
\usepackage{amssymb}
\usepackage{latexsym}
\RequirePackage{natbib}

 \usepackage{graphicx}

\usepackage{pstricks}

\usepackage[applemac]{inputenc}
\usepackage[T1]{fontenc}

\usepackage{amsfonts}
\usepackage{amsthm}
\usepackage[mathcal]{eucal}

\newtheorem{theoreme}{Theorem}[section]
\newtheorem{lemme}[theoreme]{Lemma}
\newtheorem{proposition}[theoreme]{Proposition}

\newtheorem{definition}[theoreme]{Definition}
\newtheorem{Remarque}[theoreme]{Remark}

\newtheorem{example}[theoreme]{Example}

\def\R{\mathbb R}
\def\N{\mathbb N}

\def\hop{{\noindent}}

\newcommand{\vs}[1]{\vspace{#1}}

\newcommand{\begitem}{\begin{itemize}}
\newcommand{\finit}{\end{itemize}}

\newcommand{\zun}{\vs{0.1cm}} 
 
\newcommand{\zdeux}{\vs{0.2cm}} 

\newcommand{\bx}{\mathbf{x}}

\newcommand{\by}{\mathbf{y}}

\newcommand{\bz}{\mathbf{z}}
\newcommand{\bX}{\mathbf{X}}
\newcommand{\bY}{\mathbf{Y}}
\newcommand{\bU}{\mathbf{U}}

\newcommand{\bdelta}{\boldsymbol{\delta}}

\newcommand{\CC}{\mathcal{C}}

\newcommand{\btz}{\tilde{\bz}}

\newcommand{\E}{\mathbb{E}}

\DeclareMathOperator{\support}{supp}
\DeclareMathOperator{\BC}{BC}

\begin{document}

\begin{frontmatter}

\title{Ergodic Properties of weak Asymptotic Pseudotrajectories for Set-valued Dynamical Systems }

\runtitle{WAPT relative to Set-valued Dynamical Systems}

\begin{aug}
  \author{\fnms{Mathieu}  \snm{Faure}\thanks{Institut de Math\'ematiques, Universit\'e de Neuch\^{a}tel, Rue Emile-Argand 11, 2000 Neuch\^atel, Switzerland.} \ead[label=e1]{mathieu.faure@unine.ch}}
  \and
\author{\fnms{Gregory} \snm{Roth}\thanks{Institut de Math\'ematiques, Universit\'e de Neuch\^{a}tel, Rue Emile-Argand 11, 2000 Neuch\^atel, Switzerland. Tel. $+41 32 7182828$. gregory.roth@unine.ch}  \ead[label=e2]{gregory.roth@unine.ch}}

\runauthor{M. Faure and G. Roth}
%
%
%
%

\end{aug}

\begin{abstract}
A successful method to describe the asymptotic behavior of various  deterministic and stochastic processes such as asymptotically autonomous differential equations or stochastic approximation processes is to relate them to an appropriately chosen limit semiflow. Bena\"im and Schreiber (2000) define a general class of such stochastic processes, which they call weak asymptotic pseudotrajectories and study their ergodic behavior. In particular, they prove that the weak$^*$ limit points of the empirical measures associated to such processes are almost surely invariant for the associated deterministic semiflow. Continuing a program started by Bena\"im, Hofbauer and Sorin (2005), we generalize the ergodic results mentioned above to weak asymptotic pseudotrajectories relative to set-valued dynamical systems. 
\end{abstract}



\end{frontmatter}

\section{Introduction}

Let $(E,d)$ be a complete metric space and consider a flow $\phi$ on $E$, namely $\phi : E \times \R \rightarrow E$ is continuous and satisfies the following properties:
\begitem
\item[$(i)$] for any $x \in E$, $\phi(x,0) = x$,
\item[$(ii)$] for any $t,s \in \R$ and $x \in E$, $\phi(\phi(x,t),s) = \phi(x,t+s)$.
\end{itemize}
In the sequel, we will prefer the notation $\phi_t(x)$ instead of $\phi(x,t)$. A continuous function $X : \mathbb{R}_+ \longrightarrow E$ is an asymptotic pseudotrajectory (APT) for the flow $\phi$ if
\begin{equation} \label{eq:APT} \lim_{t \to \infty} \sup_{s \in [0,T]} d(X(t+s), \phi_s(X(t))) = 0,
\end{equation}
for any $T >0$. Heuristically this means that, for any $T>0$, the curve joining $X(t)$ to $X(t+T)$ shadows the trajectory of the solution starting from $X(t)$ with arbitrary accuracy, provided $t$ is large enough. This concept has initially been introduced in \cite{Ben96} and \cite{BenHir96}, where the authors proved that the asymptotic behaviors of an APT can be described with a great deal of generality through the study of the asymptotics of the flow $\phi$. One of their main results is the characterization of the  limit set of an APT, in the sense that it is  \emph{internally chain transitive}, i.e. compact, invariant  and contains no proper attractor for the restricted flow (the terminology comes from  the notion of \emph{chain recurrence} introduced by  \cite{Con78}, see also \cite{Bow75}). 

\noindent Consequently, this result turns out to be a particularly useful tool for analyzing the long term behavior of  a large class of "perturbed" systems, whose solutions are APTs relative to some "unperturbed" flow. For instance,  given an asymptotically autonomous differential equation, its solution trajectories  are APTs relative to the flow induced by its limit autonomous differential equation (see \cite{BenHir96}). Also, under the right assumptions, the paths of a stochastic approximation process with decreasing step size are almost surely APTs for the flow induced by the mean ODE (see e.g. \cite{Ben99} or \cite{Pem07} for comprehensive overviews on the topic).

In  \cite{BenSch00}, the authors investigate the ergodic or statistical behavior of APTs for a flow. In fact, they prove their main result for a more general class of stochastic processes that they call \emph{weak asymptotic pseudotrajectories} (WAPT). Given a probability space $(\Omega, \mathcal{F}, \mathbb{P})$ and a nondecreasing family of sub-$\sigma$-algebras $(\mathcal{F}_t)_{t \ge 0}$, a process
$$\bX : \R_+ \times \Omega \longrightarrow E,$$
is called a \emph{weak asymptotic pseudotrajectory}(WAPT) for $\phi$ provided that
\begin{enumerate}
 \item [(i)] It is progressively measurable : $\bX_{|[0,T] \times \Omega}$ is $\mathcal{B}_{[0,T]} \times \mathcal{F}_T$ measurable for all $T>0$,
\item [(ii)] for each $\alpha >0$ and $T>0$,
$$\lim_{t \to \infty}\mathbb{P}(\sup_{s \in [0,T]}d(\bX(t+ s),\phi_s(\bX(t))) \ge \alpha | \mathcal{F}_t) =0$$
almost surely.
\end{enumerate}
 Notice that, if a random process $\bX$ is almost surely an APT, namely for almost every $\omega \in \Omega$, (\ref{eq:APT}) holds, then $\bX$ is  a weak asymptotic pseudotrajectory. One should be aware that the characterization of limit sets no longer holds for a WAPT. By this we mean that they are not almost surely internally chain transitive in general.
One of the main result of the paper quoted above is the following
\zun

\hop {\bfseries Theorem A} (\cite{BenSch00}, Theorem 1) {\itshape Given a WAPT $\bX$ for a flow $\phi$, the weak$^*$ limit points of the empirical measures 
$$ \mu_t(\omega) := \frac{1}{t}\int_0^t\delta_{\bX(s, \omega)}ds$$ (where $\delta_p$ is a Dirac measure at the point $p$) are almost surely invariant measures for $\phi$. In particular, if the family of measures $\{ \frac{1}{t}\int_0^t\delta_{\bX(s)}ds\}_{t \ge 0}$ is tight, this implies that, with probability one, the process $\bX(t)$ spends most of its time in any arbitrarily small neighborhood of the Birkhoff center of $\phi$.
}

\zdeux

 In \cite{BHS1}, the authors generalized the notion of asymptotic pseudotrajectory to  set-valued dynamical systems $\Phi$ (induced for instance by a differential inclusion, see Section 2 for a general definition of set-valued dynamical systems and section 4 for the particular case of differential inclusions) and extended  the characterization of limit sets. In this paper, we generalize the notion of WAPT to set-valued dynamical systems and extend Theorem A to these more general settings. 

The paper is organized as follows. The first step (see Section 2) consists in defining properly the notions of invariant measure and Birkhoff center for set-valued dynamical systems. To this end, we heavily rely on \cite{AkiMil99}. In this paper, the authors prove the equivalence between various definitions of an invariant measure for discrete time set-valued dynamical systems induced by closed relations. We give two equivalent definitions in the continuous time case (see Theorem \ref{th:eqdef}). In the same framework, \cite{AubFra91} prove a  Poincar\'e recurrence Theorem. We give a topological version of this result in the continuous time case (see Theorem \ref{poincaresetvalued}).
In Section 3, we prove the main result about the ergodic behavior of WAPTs   (Theorem \ref{thmwapt}). We then give some examples of WAPTs in Section 4, in particular stochastic approximation algorithms relative to a differential inclusion. The proofs of several technical results are postponed to the appendix to ease the reading.

\section{Ergodic theory for set-valued dynamical systems}

Set-valued dynamical systems, which are often referred to as
general control systems or general dynamical systems, are  generally used to describe multi-valued differential equations, including differential inclusions (see \cite{BiaZec81}, \cite{LiZha02}, \cite{BHS1,BHS2} or \cite{Nie09}) and control systems (see \cite{Rox65} or \cite{Klo75,Klo78}). The literature on the subject is abundant and the terminology sometimes differs among authors. 

\subsection{Generalities}

Let us first recall some classical notions.  The \emph{Hausdorff distance} between two nonempty closed sets $A$ and $B$ in $E$ is given by
\[D_H(A,B) := \max \left\{d_H(A,B),d_H(B,A) \right\},\]
where $d_H$ is the \emph{Hausdorff semidistance}:
\[d_H(A,B) := \sup_{a \in A} d(a,B).\]

 Let $\mathcal{C}(\mathbb{R},E)$ denote the space of continuous $E$-valued applications, endowed with the topology of uniform convergence on compacts. This topological space is metrizable with the distance $D$, given by
\[D(\bx,\by) := \sum_{k \in \mathbb{N}} \frac{1}{2^k} \min \{1, \sup_{t \in [-k,k]} d(\bx(t),\by(t)) \},\]
which makes it complete (since $E$ is complete).

\begin{definition}\label{setvalueddynamic}
A set-valued map  $\Phi : \R_+ \times E \rightrightarrows E$\footnote{this means that $\Phi$ is a map from $\R_+ \times E$ to $2^E$. Set-valued maps are also called \emph{relations} in the sequel.} with nonempty and closed values is called a \emph{set-valued dynamical system}(SVDS) on $E$ provided that 

\begitem
\item[$a)$] $\forall x \in E$, $\Phi_0(x) = \{x\}$,
\item[$b)$] $\forall x \in E, \; \forall s,t \in \R^+$, $\Phi_t (\Phi_s(x)) = \Phi_{t+s}(x)$,
\item[$c)$] for any $x \in E$, $t \mapsto \phi_t(x)$ is a continuous map for the Hausdorff distance,
\item[$d)$] for any $t \in \R_+$, the map $x \mapsto \phi_t(x)$ is upper semicontinuous, i.e. for any $x_0 \in E$, for any $\varepsilon >0$ there exists $\delta >0$ such that $d_H(\Phi_t(x),\Phi_t(x_0))< \varepsilon$ for any $x \in B(x_0,\delta)$. 
\end{itemize}
\end{definition}

\begin{definition}
A function $\mathbf{z}:[0,T] \rightarrow E$ is a partial solution relative to the set-valued dynamical system $\Phi$ if it satisfies $\mathbf{z}(t) \in \Phi_{t-s}(\mathbf{z}(s))$ for $s,t \in [0,T]$, $s\le t$. The set of such solutions is called $S_{\Phi}^{[0,T]}$. We call $S_{\Phi}^{[0,T]}(A)$ the set of partial solutions on $[0,T]$, starting in $A$. We denote by $S^+_{\Phi}$ the set of half solutions $S_{\Phi}^{[0,\infty)}$.
\end{definition}

 The partial solutions are continuous functions and, given $0 \leq t$ and two points $x,y \in E$ such that $y \in \Phi_{t}(x)$, there exists at least one partial solution $\mathbf{z}$ on $[0,t]$ such that $\mathbf{z}(0) = x$ and $\mathbf{z}(t) = y$ (see \cite{Rox65}). By the above remark, for all $x\in E$, $S_{\Phi}^{[0,t]}(x) \neq \emptyset$ which implies that $S_{\Phi}^{+}(x) \ne \emptyset$. A function $\mathbf{z}: \R \rightarrow E$ which satisfies $\mathbf{z}(t) \in \Phi_{t-s}(\mathbf{z}(s)), \; \forall s,t \in \R,  s  \leq t $, is called an \emph{entire solution} of $\Phi$. $S_{\Phi}$ is the set of all entire solutions, $S_{\Phi}(A)$ the subset of entire solutions starting from $A$. We will say that the set-valued dynamical system $\Phi$ is \emph{complete} if, for any $x \in E$, there exists an entire solution with initial condition $x$, i.e., if $S_{\Phi}(x) \neq \emptyset, \, \forall x \in E$. For example, we will see in Section \ref{sec_wapt_DI} that the SVDS induced by a standard differential inclusion is complete.

\zun

\hop From now, let us assume that $E$ is compact and call it $M$ to avoid confusions. Therefore, we consider a SVDS $\Phi$ defined on $M$. 

\begin{proposition}\label{prop:sphi_nonvide}
The set of solutions $S_{\Phi}$ is a nonempty compact subset of $\mathcal{C}(\mathbb{R},M)$.
\end{proposition}
\hop The compactness is a well-known consequence of Barbashin's Theorem\footnote{Which states that, for any $t \geq 0$, the set $S_{\Phi}^{[0,t]}(A)$ is compact in $\mathcal{C}([0,t],\R)$} (see for instance \cite{Bar48} or \cite{AubCel84}). We prove the existence of at least one entire solution in the Appendix \ref{appendix:sphi_nonvide}.


\begin{Remarque} Points $c)$ and $d)$ in Definition \ref{setvalueddynamic} imply that the map $(t,x) \mapsto \Phi_t(x)$ is jointly upper semicontinuous. Since $M$ is compact, upper semicontinuity is equivalent to saying that the graph of $\Phi$:
\[Gr(\Phi) := \left\{(t,x,y) \, : \; \, y \in \Phi_t(x) \right\}\]
is closed: if $(t_n,x_n,y_n) \rightarrow (t,x,y)$ and $y_n \in \Phi_{t_n}(x_n)$ then $y \in \Phi_t(x)$. 
\end{Remarque}


For our purpose, we need to give a proper definition of an invariant measure relative to set-valued dynamical systems. Recall that, if $\theta = \{\theta_t\}_{t\in \R}$ (resp. $\theta = \{\theta_t\}_{t\in \R_+}$) is a continuous flow (resp. semiflow) on a metric space $\mathcal{X}$, a measure $\mu$ on $\mathcal{X}$ is called \emph{$\theta$-invariant} if $\mu (\theta_t^{-1}(A)) = \mu(A)$ for all Borel set $A \in \mathcal{X}$ and $t \in \R$ (resp. $t \in \R_+$). 

Let $\mathcal{X}$ and $\mathcal{X}'$ be two metric spaces and $f$ be a Borel measurable map from $\mathcal{X}$ to $\mathcal{X}'$. We denote by $\mathcal{P}(\mathcal{X})$ the set of probability measures on $\mathcal{X}$. Then we can define the map $f^*: \mathcal{P}(\mathcal{X}) \rightarrow \mathcal{P}(\mathcal{X}')$ by
\[f^*(\mu)(A') := \mu(f^{-1}(A)), \]
for any $\mu \in \mathcal{P}(\mathcal{X})$ and any Borel set $A'$ in $\mathcal{X}'$. The \emph{support} of $\mu \in \mathcal{P}(\mathcal{X})$, i.e. the smallest closed set $S \subset \mathcal{X}$ satisfying $\mu(S) = 1$, is denoted by $\support(\mu)$. In the sequel,  given  a closed subset  $S$ of $\mathcal{X}$, we will sometimes assimilate a probability measure with support in  $S$ to its restriction to the metric space $(S,d)$.

\subsection{Invariant measures for a closed relation}

Let  $F : M \rightrightarrows M$ be a closed relation with nonempty values, which means that the graph of $F$,
\[Gr(F) := \left\{ (x,y) \in M \times M \mid \; \, y \in F(x)\right\}\]
is closed. Let $M^{\mathbb{Z}}$ be the set of bi-infinite sequences in $M$. The relation $F$ induces a discrete time set-valued dynamical system on $M$, whose set of solutions is the nonempty set
\[S_F := \left\{\mathbf{x} \in M^{\mathbb{Z}}: \; \, \mathbf{x}_{i+1} \in F(\mathbf{x}_i) \; \, \forall i \in \mathbb{Z}  \right\}.\] 
In order to define invariant measures in this discrete case, we follow \cite{AkiMil99}.  The shift homeomorphism $\tilde{\Theta} : M^{\mathbb{Z}} \longrightarrow M^{\mathbb{Z}}$ is defined by $\big(\tilde{\Theta} (\mathbf{x})\big)_i = \mathbf{x}_{i+1}$, for $i\in \mathbb{Z}$. 
Notice that $M^{\mathbb{Z}}$, equipped with the product topology is metrizable via the following distance:
\[\delta(\mathbf{x},\mathbf{y}) := \sup_{k \in \mathbb{Z}} \min \left\{d(\mathbf{x}_k,\mathbf{y}_k),1/|k| \right\}\]

 The following theorem is due to Akin and Miller (see \cite{AkiMil99}, Theorem 3.2). 
\begin{theoreme}
A probability measure $\mu \in \mathcal{P}(M)$ is called an \emph{invariant measure} for $F$ if it satisfies the following equivalent conditions.
\begin{enumerate}
\item For every Borel set $A \subset M$
$$\mu(A) \le \mu(F^{-1}(A)).$$
\item There exists a Markov kernel $\kappa : M \longrightarrow \mathcal{P}(M)$ satisfying
$$x \in \support(\mu) \Rightarrow \support(\kappa(x,\cdot)) \subset F(x)$$
and
$$\mu(\cdot) = \kappa^{*}(\mu)(\cdot) := \int_M \kappa(x,\cdot) \mu(dx).$$
\item There exists  $\tilde{\mu} \in \mathcal{P}(M\times M)$ the space of probability measures on $M \times M$ such that 
$$\support(\tilde{\mu}) \subset Graph(F)$$
and
$$\mu = \pi_{1}^*(\tilde{\mu}) = \pi_{2}^*(\tilde{\mu}),$$
where $\pi_i : M \times M \longrightarrow M$ is the ith projection.
\item There exists a probability measure $\nu$ on $M^{\mathbb{Z}}$ which is invariant with respect to the shift homeomorphism $\tilde{\Theta}$, satisfying
$$\support(\nu) \subset S_{F}$$
and
$$\mu = \tilde{\pi}_0^*(\nu),$$
where $\tilde{\pi}_0: \mathbf{x} \in M^{\mathbb{Z}} \mapsto \mathbf{x}_0 \in M$. 
\end{enumerate}

The set $\mathcal{P}_F(M)$ of $F$-invariant measures is nonempty, compact and convex in $\mathcal{P}(M)$. 
In general, if $\mu \in \mathcal{P}_F(M)$ then
$$\support(\mu) \subset \tilde{\pi}_{0}^*(S_F).$$
\end{theoreme}

\subsection{Invariant measures for set-valued dynamical systems}

For a  set-valued dynamical system $\Phi$, we now give two equivalent definitions of an invariant measure. Notice that, by definition, for any $t \in \mathbb{R}_+$, the set-valued map $\Phi_t: M \rightrightarrows M$ is a closed relation with nonempty images (hence its set of solutions $S_{\Phi_t}$ is nonempty).  Let us introduce the Lipschitz map $\pi_{0} : \mathcal{C}(\R, M) \rightarrow M$ defined by $\pi_{0}(\mathbf{y}) = \mathbf{y}(0)$ and the translation flow $\Theta: \mathbb{R} \times \mathcal{C}(\mathbb{R},M) \rightarrow \mathcal{C}(\mathbb{R},M)$, which associates to a real number $t$ and  $\mathbf{y} \in \mathcal{C}(\mathbb{R},M)$ the translated map $\Theta_t(\mathbf{y})$, defined by
\[\Theta_t(\mathbf{y})(s) = \mathbf{y}(t+s).\]

\begin{theoreme} \label{th:eqdef} Given a probability measure $\mu \in \mathcal{P}(M)$, the two following statements are equivalent:
\begitem
\item[$(i)$] for any $t\ge 0$, there exists a probability measure $\nu_t$ on $S_{\Phi_{t}}$ such that
\begin{enumerate}
\item $\tilde{\pi}_{0}^{*}(\nu_t) = \mu$, 
\item $\nu_t$ is $\tilde{\Theta}$-invariant,
\end{enumerate}
\item[$(ii)$] there exists a probability measure $\nu$ on $S_{\Phi}$ such that
\begin{itemize}
\item[$a)$] $\pi_{0}^{*}(\nu)=\mu$, 
\item[$b)$] $\nu$ is $\Theta$-invariant.
\end{itemize}
\end{itemize}
Such a probability measure will be called an invariant measure for the set-valued dynamical system $\Phi$. We call $\mathcal{P}_{\Phi}(M)$ (or $\mathcal{P}_{\Phi}$) the set of invariant measures for $\Phi$ . It is a nonempty compact convex subset of $\mathcal{P}(M)$. 
\end{theoreme}


\hop {\bfseries Proof.} In order to prove $(i) \Rightarrow (ii)$, we define, for all $t>0$, a new relation 
$$K_t : Gr(\Phi_t) \rightrightarrows S_{\Phi}^{[0,t]}$$
which associates, to $(x,y) \in Gr(\Phi_t)$, the set  
$$K_t(x,y)=\{\mathbf{z} \in S_{\Phi}^{[0,t]} : \mathbf{z}(0) =x \ ; \mathbf{z}(t) =y \}.$$
Notice  that $K_t(x,y)$ is nonempty and $K_t$ is a closed relation for all $t>0$. Indeed, assume that the sequence $(x_n,y_n)_n$ converges to $(x,y) \in Gr(\Phi_t)$ and that, for all $n$, $\mathbf{z}^n \in K_t(x_n,y_n)$ and $\mathbf{z}^n \rightarrow \mathbf{z} \in S_{\Phi}^{[0,t]}$. We  easily have $\mathbf{z}(0) = x$ and $\mathbf{z}(t)=y$ and so $\mathbf{z} \in K_t(x,y)$.

 Since $ S_{\Phi}^{[0,t]}$ is compact, closedness of $K_t$ is equivalent to upper semicontinuity, which can also be stated: for any $A$ closed in $S_{\Phi}^{[0,t]}$,  $K_t^{-1}(A)$ is closed in $M \times M$. In particular, $K_t$ is measurable: for any closed $A \subset S_{\Phi}^{[0,t]}$, $K_t^{-1}(A)$ is a Borel set.  By Theorem 8.1.3 of \cite{AubFra90} we can therefore choose, for all $t>0$, a measurable selection of $K_t$, 
$$\kappa_t :  Gr(\Phi_t) \rightarrow S_{\Phi}^{[0,t]}.$$

 Let $t >0$ be fixed for now. We claim that there exists a measurable application $h_t: S_{\Phi_t} \rightarrow S_{\Phi}$ which conjugates the shift operators $\tilde{\Theta}$ and $\Theta_t$:
\begin{equation} \label{conjugacy}
h_t \circ \tilde{\Theta} = \Theta_t \circ h_t.
\end{equation}

In order to prove (\ref{conjugacy}), we now define two sets :
$$A_t := \{ (x_n,y_n)_{n \in \mathbb{Z}} \in Gr(\Phi_t)^{\mathbb{Z}}, \, \mid \; y_n = x_{n+1} \, \mbox{ for all } n \};$$
and
$$B_t := \{ (\mathbf{z}^n)_{n \in \mathbb{Z}} \in (S_{\Phi}^{[0,t]})^{\mathbb{Z}} \, \mid \; \mathbf{z}^{n+1}(0) = \mathbf{z}^n(t)\}.$$

Notice that $S_{\Phi_t}$, endowed with the metric $\delta$ (see previous section)
is topologically equivalent to $A_t$, seen as a subset of the product space $Gr(\Phi_t)^{\mathbb{Z}}$, equipped with the induced product topology. Similarly, $(S_{\Phi},D)$ is topologically equivalent to the set $B_t$, understood as a subset of the product space $(S_{\Phi}^{[0,t]})^{\mathbb{Z}}$, equipped with the product of the uniform convergence topology on $[0,t]$.
We now construct a measurable function $\tilde{h}_t$  from $A_t$ to $B_t$ (to which can be associated a measurable function $h_t$ from $S_{\Phi_t}$ to $S_{\Phi}$) the following way: let $(\mathbf{x},\mathbf{y})= (x_n,y_n)_n$ be in $A_t$. Then $\tilde{h}_t (\mathbf{x},\mathbf{y})$ is given by $\mathbf{z} = (\mathbf{z}^n)_n$, where 
\[\mathbf{z}^n(s) = \kappa_t(x_n,y_n)(s), \; \forall s \in [0,t].\]
In other terms, $\tilde{h}_t$ is the countable product of the measurable map $\kappa_t$ and therefore is measurable.

 Now the corresponding map $h_t: S_{\Phi_t} \rightarrow S_{\Phi}$ is also measurable. To understand why the conjugacy (\ref{conjugacy}) holds, let us give some insights on the map $h_t$: given  $(x_n)_n \in S_{\Phi_t}$, we consider, for each couple $(x_n,x_{n+1})$ its image by $\kappa_t$ (which is a partial solution curve of length $t$) and then build a solution $\mathbf{z} \in S_{\Phi}$ by joining together these partial solution curves. The conjugacy is a clear consequence of this construction.

 By assumption, for any $t\ge0$, there exists a probability measure $\nu_t$ on $S_{\Phi_t}$ such that
\begin{enumerate}
\item $\tilde{\pi}_{0}^{*}(\nu_t) = \mu$
\item $\nu_t$ is invariant for $\tilde{\Theta}$.
\end{enumerate}
Let $(t_n)_n$ be a strictly decreasing sequence, converging to $0$ and $\nu$ be a limit point of $(h_{t_n}^{*}(\nu_{t_n}))_n$ (such a point exists since the considered sequence is tight). Notice that $\nu$ sits on $S_{\Phi}$. We can assume, without loss of generality,  that $\nu = \lim_{n \to \infty} h_{t_n}^{*}(\nu_{t_n})$. Since $\pi_0 \circ h_t = \tilde{\pi}_0$, we have  $\pi_0^*(\nu) = \mu$. There remains to prove that $\nu$ is $\Theta$-invariant. Let $f$ be a continuous function on $S_{\Phi}$ and $T\ge0$. It is sufficient to prove
\begin{equation}\label{equationinvariant}
\int_{S_{\Phi}} f(\mathbf{z}) \nu(d\mathbf{z}) = \int_{S_{\Phi}} f(\Theta_T(\mathbf{z}))\nu(d\mathbf{z}).
\end{equation}
We call 
$$\xi_n =  \int f(\mathbf{z}) (h_{t_n}^{*}(\nu_{t_n}))(d\mathbf{z}) \ \text{and } \ \xi'_n = \int f\circ \Theta_T(\mathbf{z}) (h_{t_n}^{*}(\nu_{t_n}))(d\mathbf{z}).$$
The two members of equation (\ref{equationinvariant}) can be rewritten in the form
$$ \xi := \lim_{n \to \infty} \xi_n \
\text{ and }  \  \xi' := \lim_{n \to \infty} \xi'_n.$$
Recall that $\Theta_t \circ h_t = h_t \circ \tilde{\Theta}$ for all $t\ge0$. Call  $s_n = T - [\frac{T}{t_n}]t_n$.  Then
\begin{eqnarray*}
\Theta_T \circ h_{t_n} &=& \Theta_{s_n} \circ (\Theta_{t_n})^{[\frac{T}{t_n}]} \circ h_{t_n} \\
&=&  \Theta_{s_n} \circ   h_{t_n} \circ \tilde{\Theta}^{[\frac{T}{t_n}]}.
\end{eqnarray*}
Since $\nu_{t_n}$ is $\tilde{\Theta}$-invariant, we get :
$$\xi'_n= \int f\circ \Theta_{s_n}(\mathbf{z}) (h_{t_n}^{*}(\nu_{t_n}))(d\mathbf{z}).$$
 Now we prove that $|\xi_n - \xi'_n|$ converges to zero. Pick $\varepsilon >0$.  Since $S_{\Phi}$ is compact, $(t,\mathbf{z}) \mapsto \Theta_{t}(\mathbf{z})$ is continuous and $s_n \rightarrow 0$, there exists $N$ large enough so that
$$| f\circ \Theta_{s_n}(\mathbf{z}) - f(\mathbf{z})| < \varepsilon,$$
for all $n\ge N$ and $\mathbf{z} \in S_{\Phi}$. Since $\nu_{t_{n}}$ is a probability measure for all $n$, we get the equation (\ref{equationinvariant}).
\zdeux

 Conversely, assume that $(ii)$ holds. There exists a probability measure $\nu$ which satisfies $a)$ and $b)$. For any $t\ge0$ we define the application $g_t : S_{\Phi} \longrightarrow S_{\Phi_t}$ which associates to $\mathbf{z}$, $g_t(\mathbf{z}) = (\mathbf{z}(kt))_{k\in \mathbb{Z} }$ and another probability measure $$\nu_t := g_t^*(\nu).$$
Since $\tilde{\pi}_{0} \circ g_{t} = \pi_0$, we have $\tilde{\pi}_{0}^*(\nu_t) = \tilde{\pi}_{0}^*(g_{t}^{*}(\nu))$. Now we show that $\nu_t$ is $\tilde{\Theta}$-invariant. Since $\tilde{\Theta}\circ g_t = g_t \circ \Theta_t $, we have
\begin{eqnarray*}
\tilde{\Theta}^{*}(\nu_t) &=& \tilde{\Theta}^{*}(g_{t}^{*}(\nu)) = (\tilde{\Theta} \circ g_{t})^{*}\\ 
&=& (g_{t} \circ \Theta_t)^{*}(\nu) = g_{t}^{*}((\Theta_{t})^{*}(\nu)) \\
&=& g_{t}^{*}(\nu)\\
&=& \nu_t
\end{eqnarray*}

The set $\mathcal{P}_{\Theta_{|S_{\Phi}}}$ of $\Theta_{|S_{\Phi}}$-invariant  probability measures is a convex and compact  nonempty set  (by Krylov-Bogolubov Theorem). Since $\mathcal{P}_{\Phi}$ is the image of $\mathcal{P}_{\Theta_{|S_{\Phi}}}$ under $\pi^*_0$, we easily check that it enjoys the same properties.
$\; \; \blacksquare$

\begin{Remarque}
If $\Phi = \phi$ is a flow on $M$ then $\pi_0$ restricts to an homeomorphism of $S_{\phi}$ to $M$. Consequently, $\mu$ is $\phi$-invariant if and only if it is the image of a $\Theta$-invariant measure (with its support in $S_{\phi}$) under $\pi_0$.
\end{Remarque}

\begin{Remarque}\label{rem:semi_inv}
Another natural way to define invariant measures relative to a SVDS is to consider the set of half solutions $S_{\Phi}^+$ instead of $S_{\Phi}$. Let us define the map $\Theta^+ : \R_+ \times \mathcal{C}(\R_+,M)  \rightarrow \mathcal{C}(\R_+,M)$ as the shift semi-flow defined similarly to $\Theta$. A probability measure $\mu$ on $M$ is then said to be a \emph{semi-invariant measure} for the  set-valued dynamical system $\Phi$ if there exists a probability measure $\nu^+$ on $\CC(\R_+,M)$ such that

\begin{itemize}
\item[(i)] \hspace{.3cm}$\support(\nu^+) \subset S_{\Phi}^+$, 
\item[(ii)] \hspace{.3cm}$\nu^+$ is $\Theta^+$-invariant, and
\item[(iii)] \hspace{.3cm}$(\pi_{0}^+)^{\ast}(\nu^+)=\mu$,  where $\pi_0^+ : \mathcal{C}(\R_+,M) \rightarrow M, \; \pi_0^+(\tilde{\mathbf{y}}) = \tilde{\mathbf{y}}(0)$.  
\end{itemize}
 
Semi-invariance is, a priori, very similar to invariance (as defined in Theorem \ref{th:eqdef}). We discuss the relationship between these two definitions in Appendix \ref{appendix:eqdef}.
\end{Remarque}

%
%

\subsection{Poincar\'e recurrence theorem for set-valued dynamical systems}

As an application of the definition of an invariant measure, we shall state a topological version of the Poincar\'e recurrence theorem for set-valued dynamical systems. Roughly speaking, this theorem says that invariant measures sit on the closure of the set of  recurrent points of the dynamic. The concept of recurrent point is closely related to the notion of \emph{$\omega$-limit} set of a point $x\in M$, defined by
\[
 \omega_{\Phi}(x) := \bigcap_{t \ge 0} \overline{\Phi_{[t,\infty)}(x)}
\]
It is characterized by the following: $y \in \omega_{\Phi}(x)$ if and only if there exists $(t_n)_n \uparrow \infty, \; (\mathbf{z}^n)_n \subset S_{\Phi}(x)$ such that $\mathbf{z}^n (t_n) \rightarrow y$ (see \cite{BHS1}).
Notice that $\omega_{\Phi}(x)$ contains the limit sets $L(\mathbf{z}) = \bigcap_{t \ge 0} \overline{\mathbf{z}([t,\infty))}$ of all solutions $\mathbf{z}$ with $\mathbf{z}(0)=x$ but is in general larger than 
\[L(x) := \bigcup_{\mathbf{z} \in S_{\Phi}(x)} L(\mathbf{z}). \]
(See \cite{BHS1}) 

In the classical framework of a flow, a point is recurrent provided it belongs to its own $\omega$-limit set and the topological version of Poincar\'e recurrence theorem is stated as follows
\begin{theoreme}[Poincar\'e]\label{poincare}
Let $(\mathcal{X},d)$ be a separable metric space and $\theta = \{\theta_t\}_{t\in \R}$ a flow on $\mathcal{X}$.  Let $\mu$ be an invariant measure for $\theta$. The closure of the set of all recurrent points,
\[
 R^{\omega}_{\theta}:= \{x \in \mathcal{X} \colon x \in \omega_{\theta}(x)\},
\]
is called the \emph{Birkhoff center} of $\theta$ and denoted $\BC(\theta)$.
Then 
\[
 \mu(\BC(\theta)) = 1.
\]
\end{theoreme}

The same definition of recurrence does not seem to fit in the set-valued framework. Intuitively, a point $x$ is recurrent if there exists an entire solution  starting from $x$, whose limit set contains $x$. Thus it is more natural to define the set of recurrent points of  $\Phi$ by
\[ R_{\Phi}:= \{x \in M \, \mid \;  x \in L(x) \}.\]
Clearly we have $R_{\Phi} \subset R^{\omega}_{\Phi} := \{x \in M \, \mid \; x \in \omega_{\Phi}(x) \}$. The closure of $R_{\Phi}$ will be called the \emph{Birkhoff center} of $\Phi$ and will be noted $\BC(\Phi)$. Notice that, if $\Phi$ is actually a flow, then $L(x) = \omega(x)$ and consequently $R_
{\Phi}^{\omega}=R_{\Phi}$.
The following statement is a  Poincar\'e recurrence theorem for set-valued dynamical system.
\begin{theoreme}\label{poincaresetvalued}
Let $\mu$ be an invariant measure for $\Phi$, then
\[
\mu(\BC(\Phi))=1. 
\]
\end{theoreme}

\hop {\bfseries Proof.} Let $\mu$ be an invariant measure for $\Phi$ and $\nu$ be an invariant measure for $\Theta_{|_{S_{\Phi}}}$ such that $\pi_0^*(\nu) = \mu$. First of all, notice that 
\[
\pi_0(\BC(\Theta)) \subset \BC(\Phi). 
\]
Indeed, pick $\mathbf{z} \in R_{\Theta_{|_{S_{\Phi}}}}^{\omega}$. There exists a sequence $t_n \uparrow + \infty$ such that $\Theta_{t_n}(\mathbf{z}) \rightarrow_n \mathbf{z}$. In particular, $\pi_0(\mathbf{z}) = \mathbf{z}(0) = \lim_n \mathbf{z}(t_n)$, which means that $\pi_0(\mathbf{z}) \in L(\mathbf{z}(0))$.
Using  the last inclusion, we get
\begin{eqnarray*}
\mu(\BC(\Phi)) &\ge& \mu(\pi_0(\BC(\Theta))  ) \\
&=& \nu( \pi_0^{-1} \circ \pi_0(\BC(\Theta) ))  \\
&\geq& \nu(\BC(\Theta) ).
\end{eqnarray*}
The last quantity is equal to one by Theorem \ref{poincare}.$\; \; \blacksquare$

\section{Ergodic properties of weak asymptotic pseudotrajectories}

\subsection{Definition}

Let $(\Omega,\mathcal{F},\mathbb{P})$ be a probability space and $\{\mathcal{F}_t : t\ge 0\}$ a nondecreasing family of sub-$\sigma$-algebras.  We define a process
$$\bX : \R_+ \times \Omega \longrightarrow M,$$
to be a \emph{weak asymptotic pseudotrajectory} (WAPT) for the SVDS $\Phi$ if it enjoys the following properties
\begin{enumerate}
 \item[(i)] for almost every $\omega$, the path $\bX(\cdot,\omega)$ is uniformly continuous;
 \item [(ii)] it is progressively measurable : $\bX_{|[0,T] \times \Omega}$ is $\mathcal{B}_{[0,T]} \times \mathcal{F}_T$ measurable for all $T>0$;
\item [(iii)] for each $\alpha >0$ and $T>0$, we have
$$\lim_{t \to \infty}\mathbb{P}\left(\inf_{\mathbf{z} \in S_{\Phi}} \sup_{s \in [0,T]} d\left(\bX(t+s),\mathbf{z}(s)\right) \ge \alpha \Big\vert \mathcal{F}_t\right) =0$$
almost surely.
\end{enumerate}

\subsection{Ergodic behavior}

 Given a weak asymptotic pseudotrajectory $\bX$ for $\Phi$ and $\omega \in \Omega$, let $\mu_t(\omega)$ denote the \emph{empirical occupation  measure} of $\bX(\cdot,\omega)$:
$$\mu_t(\omega) := \frac{1}{t}\int_0^t \delta_{\bX(s,\omega)}ds.$$

\begin{Remarque}
$\mu_t(\omega)$ is defined as the unique Borel measure on $M$ such that, for all continuous function $f$ on $M$, 
\[\frac{1}{t} \int_0^t f(\bX(s,\omega)) ds = \int_{x \in M}f(x) \mu_t(\omega)(dx).\]
\end{Remarque}

 Let $\mathcal{P}_{(\bX,\omega)} \subset \mathcal{P}(M)$ denote the set of weak$^*$ limit points of $\{\mu_t(\omega)\}_{t\ge 0}$. Notice that, since $M$ is compact,  $\mathcal{P}_{(\bX,\omega)}$ is  nonempty and compact. We now state the main result of this section:

\begin{theoreme}\label{thmwapt}
 Let $\bX$ be a weak asymptotic pseudotrajectory for the set-valued dynamical system $\Phi$. Then 
$$\mathcal{P}_{(\bX,\omega)} \subset \mathcal{P}_{\Phi}, \ \text{ almost surely.}$$
In particular, 
\[\bigcup_{\mu \in \mathcal{P}_{(\bX,\omega)}} \support(\mu) \subset \BC(\Phi), \; \mbox{almost surely}.\]
\end{theoreme}

 Before proving this theorem, we state some useful lemmas. First we introduce a family of probability measures $\{\nu_t\}_{t\ge0}$ on $\mathcal{C}(\mathbb{R},M)$ related to the family $\{\mu_t\}_{t\ge0}$ : given $\omega \in \Omega$,
$$\nu_t(\omega) := \frac{1}{t}\int_0^t \delta_{\Theta_s(\bX(\cdot,\omega))}ds,$$
where we use the convention that an element $X \in \CC(\R_+,M)$ can also be seen as an element of $\CC(\R,M)$ with $X(t) = X(0)$ for all $t<0$.
\begin{lemme}\label{lemmerelcomp}
 The set $\{ \Theta_t(\bX(\cdot,\omega)) : t\ge0\}$ is almost surely relatively compact.
\end{lemme}

\hop {\bfseries Proof.} Hypothesis (i) in the definition of a WAPT guarantees that the set of functions $\{ \Theta_t(\bX(\cdot,\omega)) : t\ge0\}$ is almost surely equicontinuous. Since $\bX(\mathbb{R},\omega) \subset M$ (which is compact), we can apply  Ascoli-Arzela's theorem to conclude.
$\; \; \blacksquare$
\begin{lemme} \label{lm:nutight}
The family $\{\nu_t(\omega)\}_{t\ge 0}$ is almost surely tight.
\end{lemme}
 
\hop {\bfseries Proof.} By Lemma \ref{lemmerelcomp}, we know that $\{ \Theta_t(\bX(\cdot,\omega)) : t\ge0\}$ is almost surely relatively compact. The  tightness follows since the support of $\nu_t$ is included in $\overline{\{ \Theta_t(\bX(\cdot,\omega)) : t\ge0\}}$, for any $t \ge0$.
$\; \; \blacksquare$
\vspace{3mm}

 The last lemma is a generalization (in the continuous case) of Theorem 6.9 in \cite{Wal00}. A short proof is provided in appendix \ref{app:proof_lemme_tec} for convenience.

\begin{lemme}\label{lemmelimite}
Let $(\mathcal{X},d)$ be a compact metric space, $\theta = (\theta_t)_{t\in \R}$ be a flow on $\mathcal{X}$ and $\{\sigma_t\}_{t\ge0}$ be a collection of probability measures on $\mathcal{X}$. Consider the family $\{\nu_t\}_{t\ge0}$ of probability measures on $\mathcal{X}$, defined by
$$\nu_t = \frac{1}{t}\int_0^t \theta_s^*(\sigma_t)ds.$$
 Then any limit point $\nu$ of $\{\nu_t\}_{t\ge 0}$ is $\theta$-invariant.
\end{lemme}

\hop We are now ready to prove our main result.
\zdeux

\hop {\bfseries Proof of Theorem \ref{thmwapt}.}

For all $\omega \in \Omega$ and for all $\mu \in \mathcal{P}_{(\bX,\omega)}$, there exist $(t_j)_{j\ge 0}$ going to infinity and a probability measure $\nu$ on $\mathcal{C}(\mathbb{R},M)$ such that
\begin{enumerate}
\item $\mu_{t_j} \rightarrow \mu$ and $\nu_{t_j} \rightarrow \nu$,
\item $\nu$ is $\Theta$-invariant,
\item $\pi_0^*(\nu)=\mu$.
\end{enumerate}
The first point is a direct consequence of the definition of $\mu$ and the tightness of $\{\nu_t\}_t$ (see Lemma \ref{lm:nutight}), the second point is a consequence of Lemma \ref{lemmelimite} and the last point follows from the continuity of the map $\pi_0$ and the fact that $\pi_0^*(\nu_t) = \mu_t$, for all $t\in \R_+$.
The set of all such $\nu$ will be called $A(\omega,\mu)$: 
\[A(\omega,\mu) = \left\{\nu \in \mathcal{P}(\mathcal{C}(\R,M)): \; \, \exists t_j \uparrow \infty \; \mbox{ such that } 1., 2. \mbox{ and } 3. \mbox{ hold}  \right\}\]
Let $A(\omega) = \cup_{\mu \in \mathcal{P}_{(\bX,\omega)}} A(\omega,\mu)$. We have $A(\omega) \subset \mathcal{P}_{\Theta}(\mathcal{C}(\mathbb{R},M))$, the set of $\Theta$-invariant probability measures on $\mathcal{C}(\mathbb{R},M)$. We now exhibit a set $\tilde{\Omega} \subset \Omega$ of full measure such that for all $\omega \in \tilde{\Omega}$ and for all $\nu \in A(\omega)$, $\support(\nu) \subset S_{\Phi}$. Let $\{C_k\}$ be the family of closed neighborhoods  of $S_{\Phi}$ defined by
$$ C_k = \{ \mathbf{z} \in \mathcal{C}(\mathbb{R},M) : D(\mathbf{z}, S_{\Phi}) \le \frac{1}{k}\}.$$

It is sufficient to find, for all $k \ge 0$, a set $\Omega_k \subset \Omega$ of full measure such that for all $\omega \in \Omega_k$,
$$\bigcup_{\nu \in A(\omega)} \support(\nu) \subset C_k.$$

Let $k \in \N$ and $N \in \N$ (large) be fixed. First of all we choose $T > 0$ such that $T/N \in \N$ and $\sum_{j= T/N}^{\infty}\frac{1}{2^j} < \frac{1}{2k}$. Let  $\delta >0$ be small enough such that $\delta \sum_{j= 0}^{T/N}\frac{1}{2^j}  < \frac{1}{2k}$. With these choices of $T$ and $\delta$, we have for all $\by \in \mathcal{C}(\R,M)$ and $s\ge 0$,
\[
 \sup_{u \in [- T/N, T/N]} d\left( \Theta_s(\bX)(u),\by (u)\right) \le \delta \ \Rightarrow \ D\left(\Theta_s(\bX), \by\right) \le \frac{1}{k}.
\]
Therefore we have for all $i\ge 1$,
\begin{equation}\label{pf_3.2:implication}
\sup_{u \in [0,T]} d\left( \Theta_{(i-1)T}(\bX)(u),\mathbf{y} (u)\right) \le \delta  \ \ \Rightarrow \ \ D\left(\Theta_s(\bX), \Theta_{s- (i-1)T}(\by)\right) \le \frac{1}{k},
\end{equation}
for all $s \in [(i-1+ 1/N)T, (i-1/N)T ]$. Since $S_{\Phi}$ is invariant for $\Theta$, (\ref{pf_3.2:implication}) implies that the event

\vspace{0.1cm}

\[
 \left\{\inf_{\mathbf{z} \in S_{\Phi}} \sup_{u \in [0,T]} d\left(\Theta_{(i-1)T}(\bX)(u),\mathbf{z}(u)\right) \le \delta \right\}
\]

\vspace{0.1cm}

\hop is contained in the event

\vspace{0.1cm}

\begin{equation}\label{pf_3.2:inclusion}
\left\{ \inf_{\mathbf{z} \in S_{\Phi}}D\left(\Theta_s(\bX),\mathbf{z}\right) \le \frac{1}{k}, \ \forall s \in [(i-1+\frac{1}{N})T,(i-\frac{1}{N})T]  \right\}.
\end{equation}

\vspace{0.2cm}

\hop For $n \ge 1$ set
\[U_n = \mathbb{I}_{\{\inf_{\mathbf{z} \in S_{\Phi}} \sup_{u \in [0,T]} d(\Theta_{(i-1)T}(\bX)(u),\mathbf{z}(u)) > \delta\}}\]
and
\[ M_n = \sum_{i=1}^{n}\frac{1}{i}\big(U_i -\E(U_i | \mathcal{F}_{(i-1)T})\big).\]
Since $M_n$ is a martingale and $\sup_n \mathbb{E} (M_{n}^{2}) \le 4\sum \frac{1}{i^2}$, Doob's convergence theorem implies that $(M_n)$ converges almost surely. Hence, by Kronecker lemma,
\begin{equation}\label{limitemartingale}
\lim_{n \to \infty}\frac{1}{n}\sum_{i=1}^{n}(U_i -\E(U_i | \mathcal{F}_{(i-1)T})\big)=0
\end{equation}
almost surely. In others terms, there exists $\Omega_{k,N} \subset \Omega$ (a subset of full measure) such that for all $\omega \in \Omega_{k,N}$, $(\ref{limitemartingale})$ holds. Now pick $\omega \in \Omega_{k,N}$ and $\nu \in A(\omega)$.
Let $n_j = [t_j /T]$. Then
\begin{eqnarray*}
\nu(C_k) &\geq& \lim_{j \to \infty}\frac{1}{n_j T}\sum_{i=1}^{n_j}\int_{(i-1+1/N)T}^{(i-1/N)T} \delta_{\Theta_s(X)}(C_k)ds \\
&\geq& (1-2/N) \lim_{j \to \infty} \frac{1}{n_j }\sum_{i=1}^{n_j} 
\frac{1}{T(1-2/N)} \int_{(i-1+1/N)T}^{(i-1/N)T} \delta_{\Theta_s(X)}(C_k)ds \\
&\geq& (1-2/N) \lim_j \frac{1}{n_j} \sum_{i=1}^{n_j} \mathbb{I}_{\{\inf_{\mathbf{z} \in S_{\Phi}}D(\Theta_s(\bX),\mathbf{z}) \le \frac{1}{k}, \ \forall s \ \in \  [(i-1+1/N)T, (i-1/N)T ]\}}\\
&\geq& (1-2/N) \left(1 - \lim_j \frac{1}{n_j} \sum_{i=1}^{n_j} U_i \right) \\
&\geq& (1-2/N)\left(1- \lim_j \left( \frac{1}{n_j} \sum_{i=1}^{n_j} \big(U_i -\E(U_i | \mathcal{F}_{(i-1)T})\big) + \frac{1}{n_j T}\sum_{i=1}^{n_j}\E(U_i | \mathcal{F}_{(i-1)T})\right)\right).
\end{eqnarray*}
The first term in the last equality is equal to zero. Additionally,  by definition of a WAPT, we may assume without loss of generality that, on $\Omega_{k,N}$,
\begin{equation}\label{hypothesewapt}
 \lim_{t \to \infty}\mathbb{P}\left(\inf_{\mathbf{z} \in S_{\Phi}} \sup_{u \in [0,T]}d(\bX(t+u),\mathbf{z}(u)) \ge \delta  \Big\vert \;  \mathcal{F}_t\right) = 0.
\end{equation}

Consequently, $\nu(C_k) \geq 1-2/N$ for all $\nu \in \cup_{\omega \in \Omega_{k,N}} A(\omega)$, which means that $\nu(C_k) = 1$ for all $\nu \in \cup_{\omega \in \Omega_{k}} A(\omega)$, with $\Omega_k := \cap_N \Omega_{k,N}$. We conclude the proof  by setting $\tilde{\Omega} = \bigcap_{k}\Omega_k$.
$\; \; \blacksquare$

\subsection{A simple deterministic example}

Notice that the main result of previous section is not useless in the case of APTs. To understand why, let us consider a set-valued dynamical system on the circle $S^1$, defined by the differential inclusion $\frac{dx}{dt} \in F(x)$ with

$$F(x) = \left\{\begin{array}{ll}
   [0,1] & \text{ if } x =0 \\
    1-x & \text{ if } 0< x < 1.
 \end{array}\right.$$
The dynamics have the following portrait

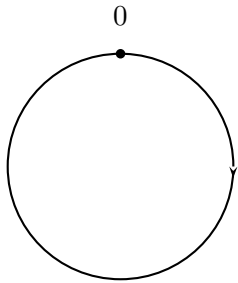
\begin{pspicture}(0,5)
\psarc{-*}(5,2.5){1.5}{0}{90}
\psarc{-<}(5,2.5){1.5}{90}{360}

\rput(5,4.5){$0$}

 \end{pspicture}
  

We easily check that the only two \emph{internally chain transitive} sets are $\{0\}$ and $S^1$. Let $X$ be a bounded APT of this dynamics. Then by Theorem 4.3 of \cite{BHS1} we know that $L(X)$ is either  $\{0\}$ or $S^1$. Moreover the only invariant measure for the dynamic is the Dirac measure on $0$. Therefore by Theorem \ref{thmwapt} we know that any limit point of the empirical measure of $X$ is $\delta_0$. In particular, this means that $X$ spends most of its time near $0$.

\section{Weak perturbed solution of a differential Inclusion} \label{sec_wapt_DI}

In this section, we give some examples of WAPTs relative to a particular case of SVDS: the set-valued dynamical systems induced by a differential inclusion. In the whole section, we are in the particular case where $E$ is  the Euclidian space $\mathbb{R}^m$.

\begin{definition}\label{standarsetvaluedmap} A set-valued map $F : \R^m \rightrightarrows \R^m$ is said to be \emph{standard} if it satisfies the following assumptions:
\begitem
\item [$(i)$] for any $x \in \R^m$, $F(x)$ is a nonempty, compact and convex subset of $\R^m$,
\item[$(ii)$] $F$ is upper semicontinuous (see Definition \ref{setvalueddynamic}),
\item [$(iii)$] there exists $c>0$ such that
\[\sup_{z \in F(x)} \|z\| \leq c (1 + \|x\|).\]
\end{itemize}
\end{definition}

 Under the above assumptions (Definition \ref{standarsetvaluedmap}), it is well known (see  \cite{AubCel84}) that the differential inclusion
\begin{equation}\label{inclusiondiff2}
\frac{d \mathbf{z}}{dt} \in F(\mathbf{z})
\end{equation}
admits at least one solution (i.e. an absolutely continuous mapping $\mathbf{z} : \R \rightarrow \R^m$ such that $\dot{\mathbf{z}}(t) \in F(\mathbf{z}(t))$ for almost every $t$) through any initial point. 
To any $x \in \mathbb{R}^m$ and $t\in \mathbb{R}_+$, we can therefore associate the nonempty set
\[ \Phi_t(x) := \left\{ \mathbf{z}(t) \mid \; \mathbf{z} \text{ is a solution of (\ref{inclusiondiff2})}, \ \mathbf{z}(0) =x  \right\}.\]
It is not hard to check that $\Phi = (\Phi_t)_{t \in \mathbb{R}}$ is  a complete set-valued dynamical system (see e.g. \cite{BHS1}). 
\zun

\begin{definition} A set $M \subset \mathbb{R}^m$ is \emph{invariant} for $\Phi$ if, for every $x \in A$, there exists an entire solution curve $\mathbf{z}$ such that $\mathbf{z}(\mathbb{R}) \subset M$. 
\end{definition}

Let $M \subset \R^m$ be a compact and invariant subset (for $\Phi$). In the following we will consider the complete  SVDS restricted to $M$, that we will also write $\Phi$ and that is defined, for all $x \in M$ and $t\in \mathbb{R}_+$, by
\[
\Phi_t(x):= \left\{ \mathbf{z}(t) \mid \; \mathbf{z} \text{ is a solution of (\ref{inclusiondiff2})}, \ \mathbf{z}(0) =x, \ \bz(\R) \subset M \right\}.
\]
Given a positive number $\delta$, let $F^{\delta}$ be the set-valued map defined by 
\begin{equation}\label{fdelta}
F^{\delta}(x) := \left\{y \mid \; \exists z \in B(x,\delta) \mbox{ such that } \; \, d(y,F(z)) < \delta \right\}.
\end{equation}

\begin{definition} Given a function $\bdelta: ]0, + \infty) \rightarrow [0,1]$ decreasing to zero as $t$ goes to infinity and a locally integrable process $\overline{\bU}: \R_+ \times \Omega \rightarrow \R^m$, we say that a process  $\mathbf{Y}: \mathbb{R}_+ \times \Omega \rightarrow M$ is a $(\bdelta,\overline{\bU})$-weak perturbed solution of the differential inclusion (\ref{inclusiondiff2}) provided  

\begin{itemize}
\item[$(i)$] $\mathbf{Y}$ is absolutely  continuous for all $\omega$,
\item[$(ii)$]   for almost every $t > 0$, 
\[\frac{d\mathbf{Y}(t)}{dt} - \overline{\bU}(t) \in F^{\bdelta(t)} (\mathbf{Y}(t)),\]
\item[$(iii)$]  for any $T> 0$ and any $\gamma >0$,
\[\lim_{t \rightarrow + \infty}  \mathbb{P}\left( \sup_{s \in [0,T]} \left\|\int_t^{t+s} \overline{\bU}(u) du \right\| \ge \gamma | \mathcal{F}_t \right)= 0,\]
almost surely.
\end{itemize}

\end{definition}

\begin{theoreme}\label{WPSWAPT}
Assume that $\bY$ is a $(\bdelta,\overline{\bU})$-weak perturbed solution  of the differential inclusion (\ref{inclusiondiff2}) and that $\overline{\bU}$ is uniformly bounded by a positive constant $C$: $\sup_{\omega \in \Omega}\sup_{t \in [0,T]} \overline{\bU}(t,\omega)\le C$. Then $\bY$ is a weak asymptotic pseudotrajectory of $\Phi$.
\end{theoreme}

\hop {\bfseries Proof.} Let $T > 0$ and define $\|F\| := \sup_{x \in M} \sup_{y \in F(x)}\|y\| < \infty$. Consider the compact set 
$$K:=  \left\{ \mathbf{y} \in Lip([0,T] , \R^m) \mid \; \, Lip(\mathbf{y}) \leq \|F\| + C + 1  \ ,   \mathbf{y}(0) \in M \right\},$$
where $Lip([0,T] , \R^m)$ denotes the set of Lipschitz  functions on $[0,T]$ and $Lip(\mathbf{y})$ is the Lipschitz constant of $\mathbf{y}$. The set $K$ is well adapted to our problem because it contains every solution curve of (\ref{inclusiondiff2}), restricted to an interval of length $T$  and every realization of any $(\bdelta,\overline{\bU})$-weak perturbed solution of the differential inclusion.

\hop For  $\delta \in [0,1]$, let us define the set-valued application (with the convention $\Lambda^0 = \Lambda$):
\begin{equation}\label{lamdadelta}
\Lambda^{\delta}:K \rightrightarrows K, \; \, \mathbf{z} \mapsto \Lambda^{\delta}(\mathbf{z}),
\end{equation} 
where $\mathbf{y} \in \Lambda^{\delta}(\mathbf{z})$ if and only if there exists an integrable function $h: [0,T] \rightarrow \R^m$ such that $h(u) \in F^{\delta}(\mathbf{z}(u)) \; \forall u \in [0,T]$ and 
\[\mathbf{y}(\tau) = \mathbf{z}(0) + \int_0^{\tau} h(u) du, \; \, \forall \tau \in [0,T].\]
Notice  that $Fix(\Lambda) := \left\{\mathbf{z} \in K \mid \; \, \mathbf{z} \in \Lambda(\mathbf{z}) \right\}$ is equal to $S_{\Phi}^{[0,T]}$, the set of partial solutions of $\Phi$ on $[0,T]$.

\hop By $(ii)$, we have 
\[\frac{d \bY(t)}{dt} - \overline{\bU}(t) \in F^{\bdelta(t)}(\bY(t)), \; \, \mbox{ for almost every } \, t>0. \]
There exists an integrable function $h: [0,T] \rightarrow \R^m$ such that $h(u) \in F^{\bdelta(t)}(\bY(t+u)) \; \forall u \in [0,T]$ and, 
 for any $\tau \in [0,T]$, 
\[\bY(t+\tau) - \int_{t}^{t + \tau} \overline{\bU}(u) du =  \bY(t) + \int_{0}^{\tau} h(u) du,\]
Hence\footnote{we call  $d_{[0,T]}$ the uniform distance on $[0,T]$}, $d_{[0,T]} \left( \bY(t + \cdot), \Lambda^{\bdelta(t)}(\bY(t+\cdot) \right) \leq \sup_{s \in [0,T]} \left\|\int_t^{t + s} \overline{\bU}(u) du  \right\|$.

Let $\alpha >0$. The following statement is an immediate consequence of Corollary 4.11 in  \cite{FauRot10}: there exists $\gamma > 0$ (which depends on $T$ and $\alpha$) and $\delta_0 > 0$ such that, for any $\delta < \delta_0$
$$d_{[0,T]}(\mathbf{z}, \Lambda^{\delta}(\mathbf{z}))< \gamma \Rightarrow  d_{[0,T]} (\mathbf{z}, S_{\Phi}) <\alpha.$$
Consequently, for $t$ large enough, 
\[d_{[0,T]} \left( \bY(t+\cdot), S_{\Phi} \right) \geq \alpha \Rightarrow d_{[0,T]} \left(\bY(t + \cdot), \Lambda^{\bdelta(t)}(\bY(t + \cdot)) \right) \geq \gamma.\]

\hop For these choices of $t$ and $\gamma$,
\begin{eqnarray*}
\mathbb{P} \left( d_{[0,T]} \left( \bY(t + \cdot), S_{\Phi} \right) \geq \alpha \mid \mathcal{F}_t \right) &\leq& \mathbb{P} \left( d_{[0,T]} \left(\bY(t+\cdot), \Lambda^{\bdelta(t)}(\bY(t+ \cdot))  \right)  \geq \gamma \mid \mathcal{F}_t \right) \\
&\leq& \mathbb{P} \left(  \sup_{s \in [0,T]} \left\|\int_t^{t+s} \overline{\bU}(u) du \right\| \geq \gamma \mid \mathcal{F}_t\right).
\end{eqnarray*}

\hop By $(iii)$, the last term tends to zero when $t$ goes to infinity and the proof is complete.$\; \; \blacksquare$

\subsection{Stochastic approximation algorithms}

Stochastic approximation algorithms were born in the early 50s through the work of \cite{RobMon51} and \cite{KieWol52}. 
Let $F : \R^m \rightarrow \R^m$ be a standard set-valued map and $M \subset \R^m$ be a compact subset invariant for the set-valued dynamical system induced by the differential inclusion (\ref{inclusiondiff2}).

\begin{definition}[Weak generalized stochastic approximation process]\label{robbinsmonro}

Let $(U_n)_n$ be an uniformly bounded $\mathbb{R}^m$-valued random process and $(F_n)_{n}$  a sequence of set-valued maps on $\R^m$. We say that $(x_n)_n$ is a \emph{generalized stochastic approximation process relative to the standard set-valued map $F$} on $M$ if the following assumptions are satisfied:

\begin{itemize} 
\item[$(i)$] we have the recursive formula
\[x_{n+1} -x_n - \gamma_{n+1} U_{n+1 }\in  \gamma_{n+1}  F_n(x_n), \]
\item[$(ii)$] the step size $(\gamma_n)_n$ is deterministic and satisfies
\[ \sum_n \gamma_n = + \infty, \; \; \lim_n \gamma_n = 0,\]
%
\item[$(iii)$] for all $n\ge 0$, $x_n \in M$,
\item[$(iv)$] for all $T>0$ and all $\gamma > 0$,
\begin{eqnarray}\label{Hypothese}
 \lim_{t \to  \infty}  \mathbb{P}\left(\sup \left\{ \left\| \sum_{i=n}^{k-1} \gamma_{i+1} U_{i+1} \right\| \mid \; \; k \, \mbox{ such that } \; \sum_{i=n}^{k-1} \gamma_i \leq T\right\}\ge \gamma \ \Big{\vert} \  \mathcal{F}_n \right)= 0. 
\end{eqnarray}
\item[$(v)$] for any $\delta > 0$, there exists $n_0 \in \mathbb{N}$
such that
\[\forall n \geq n_0, \; \, F_n(x_n) \subset F^{\delta}(x_n).\]
\end{itemize}

\end{definition}

\begin{Remarque}\label{lemmebruit}
Let $(U_n)_n$ be a $\mathbb{R}^m$-valued random process adapted to the filtration $(\mathcal{F}_n)_n$ such that
\begin{itemize}
 \item [(i)] $\E (U_{n+1} | \mathcal{F}_n) = 0$;
 \item [(ii)]for all $T > 0$, we have 
\[\lim_{R \rightarrow + \infty} \sup_n \mathbb{E}\left( \left\|  U_{n+1} \right\| \mathbf{1}_{\{\left\|  U_{n+1} \right\|\ge R\}} | \mathcal{F}_n \right) = 0.\]
\end{itemize}
Then for all $T>0$ and all $\gamma > 0$, Property (\ref{Hypothese}) is satisfied (see \cite{BenSch00}).
\end{Remarque}

Consider a weak generalized stochastic approximation process $(x_n)_n$. Set $\tau_n := \sum_{i=1}^n \gamma_i$ and $m(t) := \sup \{ j \mid \tau_j \leq t \}$. We call $\bX$ the continuous time affine interpolated process induced by $(x_n)_n$, $\overline{\gamma}$ the piecewise constant deterministic process induced by $(\gamma_n)_n$ and  $\overline{\bU}$ is the piecewise constant continuous time process associated to $(U_n)_n$:

$$\bX (\tau_i + s) = x_i + s \frac{x_{i+1} - x_i}{\gamma_{i+1}}, \mbox{ for }  \; \, s \in [0,\gamma_{i+1}] \ \overline{\gamma}(\tau_i + s) := \gamma_{i+1} \mbox{ for } \; s \in [0, \gamma_{i+1}[,$$
 
$$\mbox{and } \overline{\bU}(t) := U_{n+1}, \; \; \mbox{ for } t \in [\tau_n, \tau_{n+1}].$$

\begin{theoreme} \label{DSAperturbed} The interpolated process $\bX$ is a WAPT. Hence the conclusions of Theorem \ref{thmwapt} hold.
\end{theoreme}

\hop {\bfseries Proof.} By straightforward computations (see the proof of  proposition 1.3 in Bena\"im et al. \cite{BHS1}), it is not difficult to see that $(\bX(t))_t$ is a weak perturbed solution associated to $\overline{\bU}$ and 
\[\bdelta (t) := \inf \left\{\delta > 0 \mid \; \; \tau_n \geq t \Rightarrow \, F_n(x_n) \subset F^{\delta}(x_n) \right\} + \overline{\gamma}(t) \left(\overline{\bU}(t) + c \left(1 + \sup_{x \in M} F(x) \right) \right),\]
which converges to $0$. Consequently $ \bX$ is a WAPT relative to the SVDS induced by $F$ and the proof is complete. $\; \; \blacksquare$

\appendix
\section{Proof of Proposition \ref{prop:sphi_nonvide}}\label{appendix:sphi_nonvide}

Recall that an element $\btz \in \mathcal{C}(\R_+,M)$ can also be seen as an element of $\mathcal{C}(\R,M)$, with the convention  $\btz(t)=\btz(0)$ for all $t<0$. Since $S_{\Phi}^+ \ne \emptyset$, the existence of at least one entire solution is a direct consequence of the following Proposition. 

\begin{proposition}\label{lem_sol_entiere}
 Let $(t_n) \uparrow \infty$   be a sequence of positive real numbers converging to infinity and $\tilde{\bz} \in S_{\Phi}^+$ be a solution. Then there exists a subsequence $(t_{n_k})_{k\ge 0}$ and an entire solution $\bz \in S_{\Phi}$ such that
\[
 \lim_{k \to \infty}\Theta_{t_{n_k}}(\btz) = \bz.
\]
\end{proposition}
\hop {\bfseries Proof.} Pick some positive integer $N$. By  Barbashin Theorem, the set $S_{\Phi}^{[-N,N]}$ is compact and, for $n \in \N$ large enough,  $\Theta_{t_n}(\btz) \in S_{\Phi}^{[-N,N]}$. Therefore there exist an increasing sequence $\psi^N : \N \rightarrow \N$  and $\bz^N \in S_{\Phi}^{[-N,N]}$ such that 
\[
 \lim_{n \to \infty}\Theta_{t_{\psi^N(n)}}(\btz) = \bz^N, 
\]
By the same arguments, there exist an extraction $\psi$ and $\bz^{N+1} \in S_{\Phi}^{[-N-1,N+1]}$ such that
\[
 \lim_{n \to \infty}\Theta_{t_{\psi \circ \psi^N(n)}}(\btz) = \bz^{N+1}, 
\]
and then, in particular, $\bz^{N+1}_{\vert [-N,N]} = \bz^N$. 
We set $\psi^{N+1} := \psi \circ \psi^N$ and we iterate the process. In this way we construct an entire solution $\bz \in S_{\Phi}$ such that $\bz_{\vert [-N,N]} = \bz^N \; \forall N$.

Let $(\delta_k) \downarrow 0$ be a decreasing sequence of positive real numbers. There exists a natural number $m_1$ such that for all $m \ge m_1$, 
\[
 \sup_{s \in [-1,1]}d(\Theta_{t_{\psi^1(m)}}(\btz)(s), \bz(s)) < \delta_1, \
\]
We set $n_1 := \psi^1(m_1)$. Now we define $n_k$ by induction. Fix $k >1$. There exists a natural number $m_k \ge m_{k-1}$ such that for all $m \ge m_k$, 
\[
 \sup_{s \in [-k,k]}d(\Theta_{t_{\psi^k(m)}}(\btz)(s), \bz(s)) < \delta_k, \
\]
We set $ n_k := \psi^k(m_k)$.
Therefore, by construction, we have $\lim_{k \to \infty} \Theta_{t_{n_k}}(\btz) = \bz$.
$\; \; \blacksquare$
\zdeux

\section{Proof of Lemma \ref{lemmelimite}}\label{app:proof_lemme_tec}

Let $\mu$ be a limit point of $\{\mu_t\}_{t\ge 0}$: there exists a sequence $(t_n)_n \uparrow + \infty$  such that $\mu = \lim_{n \to \infty}\mu_{t_n}$. We have to prove that, for any $T$ and any continuous function $f$ ,
\begin{equation}\label{equation1}
\int_{\mathcal{X}} f(x)d\mu(x) = \int_{\mathcal{X}} f(\theta_T(x))d\mu(x).
\end{equation}

 We have
$$\int_{\mathcal{X}} f(x)d\mu(x) = \lim_n \frac{1}{t_n}\int_{0}^{t_{n}}  \int_{\mathcal{X}} f(\theta_s(x))d\sigma_{t_{n}}(x) ds$$
and
$$\int_{\mathcal{X}} f(\theta_T(x))d\mu(x)= \lim_n \frac{1}{t_n}\int_{0}^{t_{n}}  \int_{\mathcal{X}} f(\theta_{s+T}(x))d\sigma_{t_{n}}(x) ds.$$
By Fubini's Theorem we can exchange the integral operators in both expressions. Consequently, 
\begin{eqnarray*}
&& \frac{1}{t_n} \left|   \int_{0}^{t_{n}} \int_X f(\theta_s(x))ds\big)d\sigma_{t_{n}}(x)  -  \int_{0}^{t_{n}}  \int_X  f(\theta_{s+T}(x))ds \big) d\sigma_{t_{n}}(x)\right| \\
&\le&  \frac{1}{t_n} \left| \int_{X}  \int_{0}^{t_{n}}  f(\theta_s(x))ds\big)d\sigma_{t_{n}}(x)  - \int_{X} \int_{0}^{t_{n}}  f(\theta_{s+T}(x))ds \big) d\sigma_{t_{n}}(x)\right| \\
&\le&    \frac{1}{t_n} \int_{X} \left( \left| \int_{0}^{t_{n}}  f(\theta_s(x))ds \right| + \left|\int_{T}^{t_{n}+T}  f(\theta_{s}(x)) ds \right| \right)d\sigma_{t_{n}}(x)\\
&\le&\frac{2 T ||f||_{\infty}}{t_n}.
\end{eqnarray*}
Finally, taking the limit as $n$ goes to infinity, we obtain  (\ref{equation1}).
$\; \; \blacksquare$

\section{Some remarks on semi-invariance}\label{appendix:eqdef}

In this section we show (Proposition \ref{prop:eqdefinv}) that every invariant measure for a SVDS $\Phi$ is a semi-invariant measure for $\Phi$. To do this, we start by proving two technical lemmas.  Let us define the projection $\by^+$ on $\mathcal{C}(\R_+,M) $ of an element $\by \in \mathcal{C}(\R,M) $ by
\[
\by^+(t) = \by(t), \ \forall t \ge 0,
\]
and, in the same way, $A^+ := \{ \by^+ : \by \in A\}$ for a subset $A \subset \mathcal{C}(\R,M)$. Be aware that the set $(S_{\Phi})^+$ is contained in the set $S_{\Phi}^+$ of partial solutions, but not equal in general.

\begin{proposition} \label{prop:eqdefinv}
An invariant measure $\mu$ on $M$ for the set-valued dynamical system $\Phi$ is a semi-invariant measure for $\Phi$.
\end{proposition}

\hop {\bfseries Proof.} There exists a probability measure $\nu$ on $S_{\Phi}$ such that
\begin{itemize}
\item[$a)$] $\pi_{0}^{*}(\nu)=\mu$, 
\item[$b)$] $\nu$ is $\Theta$-invariant.
\end{itemize}
We need to construct a probability measure $\nu^+$ on $\CC(\R_+,M)$ which satisfies conditions (i), (ii) and (iii) of the definition of a semi-invariant measure given in Remark \ref{rem:semi_inv}. A natural way to do this is to define $\nu^+$  as
\[
 \nu^+(A) := \nu( B_A), \text{ for all Borel sets } A \text{ of } \CC(\R_+,M),
\]
where $B_A := \{\bz \in S_{\Phi} : \bz^+ \in A\}$.
First,  we have 
\begin{eqnarray*}
\nu^+(S_{\Phi}^+) &=& \nu\left(\{\bz \in S_{\Phi} : \bz^+ \in S_{\Phi}^+\} \right)\\
&=& \nu\left(S_{\Phi} \right) = 1.
\end{eqnarray*}
which gives condition (i). Let $T>0$ and $A \subset S_{\Phi}^+$ a Borel set. Since, for any $\bz \in S_{\Phi}$, we have $\Theta^+_T(\bz^+) = (\Theta_T(\bz))^+$,
\begin{eqnarray*}
 B_{(\Theta^+_T)^{-1}(A)}&=&\{\bz \in S_{\Phi} : \bz^+ \in (\Theta^+_T)^{-1}(A)\} \\
&=& \{\bz \in S_{\Phi} : (\Theta_T(\bz))^+ \in A\} \\
&=& \Theta_T^{-1}(B_A).
\end{eqnarray*}
Therefore, as $\nu$ sits on $S_{\Phi}$ and is $\Theta$-invariant, we have
\begin{equation*}
\nu^+\left((\Theta^+_T)^{-1}(A)\right) = \nu\left( B_{(\Theta^+_T)^{-1}(A)} \right) = \nu\left( \Theta_T^{-1}(B_A) \right) = \nu\left(B_A \right)
= \nu^+(A),
\end{equation*}
which gives condition (ii). 

Let $D$ be a Borel subset of $M$. Notice that $B_{(\pi_0^+)^{-1}(D)} = \pi_0^{-1}(D) \cap S_{\Phi} = S_{\Phi}(D)$. Consequently
\[
 (\pi_0^+)^{*}(\nu^+)(D) = \nu^+((\pi_0^+)^{-1}(D)) = \nu(S_{\Phi}(D)) = \nu(\pi_0^{-1}(D)) = \mu(D)
\]
and the result holds. $\; \; \blacksquare$

Whether the converse statement is also true is an open question.

\paragraph{{\bfseries Acknowledgements}}  The authors thank Michel Bena\"im for useful advices and discussions. They also acknowledge financial support from the Swiss National Science Foundation Grant 200020 130574.

\bibliographystyle{imsart-nameyear}
\bibliography{biblio-generale}







\end{document}